\documentclass[a4paper,oneside,10pt,final,reqno]{amsart}

\usepackage{amssymb}
\usepackage[mathcal]{euscript}
\usepackage[frame,cmtip,curve,arrow,matrix,line,graph]{xy}
\usepackage{enumerate}

\usepackage[dvipdfmx]{hyperref}

\usepackage[all, warning]{onlyamsmath}

\numberwithin{equation}{section}
\allowdisplaybreaks

\theoremstyle{definition}
\newtheorem{thm}{Theorem}[section]
\newtheorem{prop}[thm]{Proposition}
\newtheorem{dfn}[thm]{Definition}
\newtheorem{fct}[thm]{Fact}
\newtheorem*{dfn*}{Definition}
\newtheorem*{rmk*}{Remark}

\newcommand{\wh}{\widehat}
\newcommand{\seteq}{\, := \, }
\newcommand{\dsimeq}{\, \simeq \, }

\newcommand{\bbC}{\mathord{\mathbb{C}}}
\newcommand{\bbG}{\mathord{\mathbb{G}}}
\newcommand{\bbZ}{\mathord{\mathbb{Z}}}

\newcommand{\shA}{\mathord{\mathcal{A}}}
\newcommand{\shD}{\mathord{\mathcal{D}}}
\newcommand{\shG}{\mathord{\mathcal{G}}}
\newcommand{\shK}{\mathord{\mathcal{K}}}
\newcommand{\shL}{\mathord{\mathcal{L}}}
\newcommand{\shM}{\mathord{\mathcal{M}}}
\newcommand{\shO}{\mathord{\mathcal{O}}}
\newcommand{\shP}{\mathord{\mathcal{P}}}
\newcommand{\shV}{\mathord{\mathcal{V}}}

\newcommand{\shHom}{\mathop{\mathcal{H}\mathit{om}} \nolimits}

\newcommand{\catC}{\mathord{\mathcal{C}}}
\newcommand{\catM}{\mathord{\mathcal{M}}}
\newcommand{\catS}{\mathord{\mathcal{S}}}
\newcommand{\catSch}{\mathord{\mathcal{S}\mathit{ch}}}

\newcommand{\frkM}{\mathord{\mathfrak{M}}}

\newcommand{\g}{\mathord{\mathfrak{g}}}

\newcommand{\Gr}{\mathord{\mathcal{G}\mathit{r}}}

\newcommand{\ch}{\text{ch}}

\newcommand{\Wch}{\mathop{\mathcal{W}^{\text{ch}}} \nolimits}
\newcommand{\Cl}{\mathop{\mathcal{C}\mathit{l}} \nolimits}
\newcommand{\Clch}{\mathop{\mathcal{C}\mathit{l}^{\text{ch}}} \nolimits}

\newcommand{\C}{\mathop{\mathrm{C}}\nolimits}
\newcommand{\GL}{\mathop{\mathrm{GL}}\nolimits}
\newcommand{\SC}{\mathop{\mathrm{SC}}\nolimits}

\newcommand{\Hom}{\mathop{\mathrm{Hom}}\nolimits}
\newcommand{\Lie}{\mathop{\mathrm{Lie}}\nolimits}
\newcommand{\Pic}{\mathop{\mathrm{Pic}}\nolimits}
\newcommand{\Spf}{\mathop{\mathrm{Spf}}\nolimits}
\newcommand{\Spec}{\mathop{\mathrm{Spec}}\nolimits}

\newcommand{\LieD}{\mathop{\mathcal{L}\mathit{ie}}\nolimits}

\newcommand{\Wedge}{{\textstyle{\bigwedge}}}

\newcommand{\inj}{\hookrightarrow}
\newcommand{\surj}{\twoheadrightarrow}
\newcommand{\simto}{\xrightarrow{\sim}}
\newcommand{\longto}{\longrightarrow}
\newcommand{\longsimto}{\xrightarrow{\ \sim \ }}
\newcommand{\longinj}{\lhook\joinrel\longrightarrow}

\newcommand{\pair}{\left< \, , \, \right>}
\newcommand{\ket}[1]{\left| #1 \right>}

\begin{document}


\title{Boson-fermion correspondence from factorization spaces}

\author{Shintarou Yanagida}
\address{Graduate School of Mathematics, Nagoya University 
Furocho, Chikusaku, Nagoya, Japan, 464-8602.}
\email{yanagida@math.nagoya-u.ac.jp}

\date{November 18, 2016}

\begin{abstract}
We give a proof of the boson-fermion correspondence
(an isomorphism of lattice and fermion vertex algebras)
in terms of isomorphism of factorization spaces.
\end{abstract}

\maketitle

\setcounter{section}{-1}
\section{Introduction}
\label{sect:intro}

\subsection{}

The notion of factorization space 
is a non-linear analog of the factorization algebra,
which was introduced by Beilinson and Drinfeld \cite{BD}
in their theory of chiral algebras developing a geometric framework 
of vertex algebras
and the study of geometric Langlands correspondence \cite{BD2}.
A factorization space over a scheme $X$ 
is a family of ind-schemes $\shG_I$ over $X^I$ for each finite set $I$
with ``factorization structure" given by isomorphisms 
between these $\shG_I$'s.
Given a factorization space on $X$,
one can apply a linearization procedure,
and if $X$ is a smooth curve
then one obtains a factorization algebra,
an equivalent notion of chiral algebra.

Although the formulation of factorization space 
looks highly complicated at first glance,
it fits moduli problems over algebraic curves very well.
It captures the intimate connection between 
two-dimensional quantum conformal field theories
and moduli problems on curves.

Along this line, one can ask a question:
are there enhanced versions of properties of vertex algebras
in the level of factorization spaces?
In this note, we consider this kind of problem 
for the boson-fermion correspondence.

Recall that the two-dimensional boson-fermion correspondence 
is stated as an isomorphism
\begin{equation}\label{eq:V=Lambda}
 V_{\bbZ} \dsimeq \Wedge
\end{equation}
between the lattice vertex algebra $V_{\bbZ}$ attached to lattice $\bbZ$ 
and the free fermion vertex super algebra $\Wedge$.
Here we borrowed the notations in \cite{FB}.
The essence of the correspondence is that 
the vertex operators
\[
 V_{\pm}(z) \seteq :e^{\pm \varphi(z)}:, \quad
 \varphi(z) \seteq q + a_0 \log(z) 
                     - \sum_{n \in \bbZ\setminus\{0\}} \dfrac{a_n}{n} z^{-n}
\]
acting on $V_{\bbZ}$
obey the same commutation relation 
with the fermionic operators $\psi(z),\psi^*(z)$
acting on $\Wedge$.
Here $a_n$'s and $q$ denote the Heisenberg generators 
with the commutation relation $[a_m,a_n]=m \delta_{m+n,0}$
and $[q,a_n]=\delta_[n,0]$.
Thus schematically we have 
\begin{align}\label{eq:Vp}
 V_+(z) \longleftrightarrow \psi(z),\quad
 V_-(z) \longleftrightarrow \psi^*(z)
\end{align}
of vertex operators.
See \cite[\S5.3]{FB} for the detail.

Our main statement in this note is Theorem \ref{thm:main}
giving an isomorphism 
\begin{equation}\label{eq:main}
 \shG(X,\bbZ) \simeq \Gr(X,\SC(2))
\end{equation}
between factorization spaces over a smooth curve $X$.
Here $\shG(X,\bbZ)$ is the factorization space 
arising from 
the Picard scheme (moduli space of line bundles) on $X$.
$\Gr(X,\SC(2))$ is the Beilinson-Drinfeld Grassmannian for 
the Clifford group $\SC(2)$.

These factorization spaces are equipped with 
factorization super linear bundles respectively,
and the twisted linearization procedure 
give the lattice and Clifford chiral algebras.
These chiral algebras correspond to the vertex algebras 
$V_{\bbZ}$ and $\Wedge$ respectively.
Thus the isomorphism \eqref{eq:main} can be considered as 
an enhanced version of 
the boson-fermion correspondence \eqref{eq:V=Lambda}.

\subsection{Organization}

Let us explain the organization of this article.

In \S\ref{sect:fs} we recall the definition of factorization spaces,
and factorization line bundles.
We also recall the Beilinson-Drinfeld Grassmann $\Gr(X,G)$ 
for a smooth algebraic curve $X$ and a reductive group $G$,
the standard example of factorization space,

\S\ref{sect:linear} reviews the linearization procedure 
which produces a chiral algebra from factorization space.

In \S\ref{sect:lattice} we recall the factorization space 
$\shG(X,\bbZ)$ arises from the Picard functor.
By the linearization procedure 
it gives the lattice vertex algebra $V_{\bbZ}$.

In \S\ref{sect:cl} we consider the factorization space 
$\Gr(X,\SC(Q))$ for the Clifford group $\SC(Q)$ 
attached to the non-degenerate quadratic form $Q$.
We explain that by the linearization procedure 
it yields the Clifford chiral algebra.

The final \S\ref{sect:main} gives the proof of the main theorem.
The proof is not difficult,
and the key observation is that 
we have a natural identification of moduli stacks 
\[
 \catM(X,\SC(Q)) \dsimeq \Pic(X) \times \Pic(X)
\]
of the moduli space of $\SC(Q)$-bundles on $X$ 
with two pieces of Picard varieties.
It reflects the correspondence \eqref{eq:Vp}
of vertex operators 
$(V_+,V_-) \leftrightarrow (\psi,psi^*)$.

The ingredients in \S\S\ref{sect:fs}--\ref{sect:lattice} are 
basically known facts, 
and the main sources of our presentation are \cite{BD,FB}.
The content in \S\ref{sect:cl} is also a variant of 
\cite[Chap.\ 20]{FB}.

\subsection{Notation}

We will work over the field $\bbC$ of complex numbers
unless otherwise stated.
The symbol $\otimes$ denotes the tensor product of 
$\bbC$-linear spaces.

For a scheme or an algebraic stack $Z$, 
$\shO_Z$, $\Theta_Z$, $\Omega_Z$ and $\shD_Z$ 
denote the structure sheaf, the tangent sheaf, the sheaf of $1$-forms 
and the sheaf of differential operators on $Z$ respectively
(if they are defined).
By ``an $\shO$-module on Z" we mean a quasi-coherent sheaf on $Z$.
By ``a $\shD$-module on Z" we mean a sheaf of \emph{right} $\shD_Z$-modules 
quasi-coherent as $\shO_Z$-modules.

For a morphism $f: Z_1 \to Z_2$, the symbols $f^{\cdot}$ and $f_{\cdot}$ 
denote the inverse and direct image functors of $\shO$-modules 
respectively.


Finally we will use the symbols 
$\shO := k[[t]]$ for the algebra of formal series
and $\shK := k((t))$ for the field of formal Laurent series.

\section{Factorization space}
\label{sect:fs}

We follow \cite[\S3.10.16]{BD}, \cite[\S 20.4.1]{FB}, \cite[Chap.\ 5]{G} 
and \cite{KV}.

\subsection{The category $\catS$ of finite sets and surjections}

We will repeatedly use the following category of sets.

\begin{dfn}\label{dfn:catS}
Let $\catS$ be the category  of finite non-empty sets and surjections.
For $\pi:J \surj I$ in $\catS$ 
and $i \in I$ we set $J_i :=\pi^{-1}(i) \subset J$.
\end{dfn}

Let $X$ be a smooth complex curve as before.
For $\pi:J \surj I$ in $\catS$, denote by
\begin{equation*} 
 \Delta^{(\pi)} \equiv \Delta^{(J/I)}: X^I \longinj X^J
\end{equation*}
an embedding of the locus such that 
$x_j=x_{j'}$ if $\pi(j) = \pi(j')$ for $j,j' \in J$.
Also set 
\[
 U^{(\pi)} \equiv U^{(J/I)} \seteq 
 \left\{(x_j)_{j \in J} \in X^J \mid 
   x_{j} \neq x_{j'} \text{ if } \pi(j) \neq \pi(j')\right\}.
\]
One can also see that $U^{(\pi)}$ is the complement of the diagonals 
that are transversal to 
$\Delta^{(\pi)}: X^I \inj X^J$. 
Denote the open embedding $U^{(\pi)} \subset X^J$ by
\[
 j^{(\pi)} \equiv j^{(J/I)} :  U^{(\pi)} \longinj X^J.
\]

\subsection{Recollection on ind-schemes}
\label{subsec:fs:ind}

Before introducing the factorization spaces,
we need to recall the notion of ind-schemes,
since otherwise we have no good examples.
We follow \cite[\S1.1]{KV} for the presentation.

For any category $\catC$,
an ind-objects of $\catC$ is a filtering inductive system 
over $\catC$.
Ind-objects form a category where a morphism 
is a collection of morphisms between the objects in the inductive systems
satisfying some compatibility conditions.
An ind-object will be represented by the symbol 
$``\varinjlim_i" C_i$ with $C_i \in \catC$.

Denote by $\catSch$ the category of separated schemes over $\bbC$.
An \emph{ind-scheme} is an ind-object of $\catSch$ represented by 
an inductive system of schemes.
Formal schemes are ind-schemes, like 
\[
 \Spf k[[t]] = ``\varinjlim_n" \Spec k[t]/(t^{n+1}).
\]
A \emph{strict ind-scheme} is an ind-scheme with an inductive system
given by closed embeddings of quasi-compact schemes.

Recall that a scheme $S$ is equivalent to the functor of points 
$F_s: T \mapsto \Hom_{\catSch}(T,S)$.
This functor is a sheaf of sets on $\catSch$ 
which is seen as the Zariski site.
Calling such a sheaf a $\bbC$-space.
An ind-scheme is a $\bbC$-spaces 
represented by an inductive system of schemes.

One can define an ind-scheme over a scheme $Z$ similarly
by replacing the category $\catSch$ 
by the category $\catSch_Z$ of schemes over $Z$.
For a morphism $f:Z_1 \to Z_2$ of schemes,
denote by $f_*$ and $f^*$ the push-forward and the pull-back of 
ind-schemes over $Z_i$'s.

\subsection{Definition of factorization space}
\label{subsec:fs:dfn}


Now we turn to our main object, the factorization space.
It is a non-linear analog of 
factorization algebras,
and may be considered as sheaves on the Ran space 
\cite{R} (see also \cite[Chap.\ 4]{BD}).
Recall once again the category $\catS$ in Definition \ref{dfn:catS}.

\begin{dfn*} 
Let $X$ be a scheme.
A \emph{factorization space} $\shG$ on $X$ consists of the following data.
\begin{itemize}
\item 
A formally smooth ind-scheme $\shG_I$ over $X^I$ for each $I \in \catS$.

\item
An isomorphism
\[
 \nu^{(\pi)} \equiv \nu^{(J/I)}: \,
 \Delta^{(\pi) *} \shG_J \longsimto \shG_I
\]
of ind-schemes over $X^I$
for each $\pi: J \surj I$ in $\catS$.

\item
An isomorphism
\[
 \kappa^{(\pi)} \equiv \kappa^{(J/I)}: \, 
 j^{(\pi) *}\bigl(\prod_{i\in I} \shG_{J_i}\bigr) 
 \longsimto j^{(\pi) *}\shG_J
\]
of ind-schemes over $U^{(\pi)}$
for each $\pi: J \surj I$ in $\catS$,
called the \emph{factorization isomorphism}.
\end{itemize}
These should satisfy the following compatibility conditions.
\begin{enumerate}
\item 
$\nu^{(\pi)}$'s are compatible with compositions of surjections $\pi$.

\item  
For any $\pi:J \surj I$ and $\rho:K \surj J$, we should have
\[
 \kappa^{(K/J)} = 
 \kappa^{(K/I)}\bigl(\boxtimes_{i \in I}\kappa^{(K_i/J_i)}\bigr).
\]
\item
For any $J \surj I$ and $K \surj J$,
we should have
\[
 \nu^{(K/J)} \Delta^{(K/J) *}\bigl(\kappa^{(K/I)}\bigr)
 = \kappa^{(J/I)} \bigl(\boxtimes_{i \in I} \nu^{(K_i/J_i)}\bigr).
\] 
\end{enumerate}
\end{dfn*}

A factorization space $\shG$ is attached with the structure morphism
$r^{(I)}: \shG_I \to X^I$ of ind-schemes.


Let us also introduce the object corresponding to 
the unit of factorization algebra.

\begin{dfn*}
A \emph{unit} of a factorization space $\shG$ on a scheme $X$ is 
the data of morphisms 
\[
 u^{(I)}: X^I \longto \shG_I
\] of ind-schemes for each $I \in \catS$ 
such that for any morphism $f: U \to \shG_{\{1\}}$ with open $U \subset X$,
$u^{(\{1\})} \boxtimes f$, 
which can be seen as a morphism 
$U^2\setminus \Delta \to \shG_{\{1,2\}}$ by $\kappa^{(\pi:\{1,2\} \surj \{1\})}$,
extends to a morphism $U^2\to \shG_{\{1,2\}}$,
and $\Delta^{(\pi),*}(u^{(\{1\})} \boxtimes f) = f$.
\end{dfn*} 


\begin{fct}\label{fct:conn}
A factorization space $\shG$ on $X$ together with a unit 
has a connection along $X$.
\end{fct}

Let us explain precisely what a connection on $\shG$ along $X$ is.
Assume that we are given
\begin{itemize}
\item 
$I$ an local Artin scheme $I$ of length $1$ 
\item
a morphism $f: I \times S \to X$
of schemes with $S$ a scheme,
\item
$g_0:I_0 \times S \to \shG$ of ind-schemes
with $I_0 := I_{\text{red}} \simeq \Spec(k)$ 
the reduced scheme of $I$
\end{itemize}
such that $r^{(1)} \circ g_0: I_0 \times S \to \shG \to X$ coincides with 
the composition $I_0 \times S \to I \times S \to X$.
A connection on $\shG$ along $X$ is equivalent to the property 
that for given $(I,f,g_0)$  there is a map $g:I \times S \to \shG$
such that $r \circ g = f$ extending $g_0$.

\begin{dfn*}
A factorization space with a unit will be called a \emph{factorization monoid}.
\end{dfn*}

We followed \cite{KV} for the terminology.
In \cite[\S3.10.16]{BD} factorization monoid is called chiral monoid.

We also need a line bundle living over factorization space.

\begin{dfn}\label{dfn:line}
Let $\shG$ be a factorization space over $X$.
A \emph{factorization line bundle} $\shL$ over $\shG$ 
is a collection of line bundles $\shL_I$ on $\shG_I$
together with isomorphisms
\[
 j^{(J/I) \, \cdot}\shL_I \longsimto 
 j^{(J/I) \, \cdot}\bigl(\otimes_{i \in I} \shL_{J_i}\bigr)
\]
over $U^{(J/I)}$
which should satisfy the factorization property.
\end{dfn}

For later use we will name a standard example of factorization space,
namely the Beilinson-Drinfeld Grassmannian $\Gr(X,G)$.
Here $G$ is a reductive group and $X$ is a smooth algebraic curve.
%
%
%
As a preparation 
we recall some facts on the moduli space of $G$-bundles on an algebraic curve.

\subsection{Moduli space of $G$-bundles on curve}
\label{subsec:Bun_G}

Let $G$ be a reductive algebraic group over $\bbC$.
Recall that the \emph{affine Grassmannian} 
\[
 G(\shK)/G(\shO) = G\bigl(k((z))\bigr)/G\bigl(k[[z]]\bigr)
\]
can be considered as the moduli space of $G$-bundles 
on the disc $D=\Spec\shO$ together with a trivialization 
on $D^{\times} = \Spec\shK$.
Strictly speaking,
$G(\shK)$ is an ind-scheme,
and $G(\shK)/G(\shO)$ is a formally smooth strict ind-scheme.

Denote by $\catM(X,G)$ the category of $G$-torsors on $X$.
By the result of \cite{DS}, we have

\begin{fct} 
\label{fct:DS}
Let $X$ be a smooth algebraic curve and $x \in X$ be a point.
\begin{enumerate}
\item 
A choice of local coordinate $z$ at $x$ 
gives an identification
\begin{align}
\label{eq:BDG:Gr-XGx}
\Gr(X,G)_{x} \seteq 
\{(\shP,\varphi) \mid 
   \shP \in \catM(X,G), \, 
   \varphi: \text{trivialization of }
      \left.\shP\right|_{X\setminus \{x\}} \} 
\longsimto
  G(\shK)/G(\shO). 
\end{align}
\item
If $G$ is semi-simple, then
any $G$-bundle on $X\setminus\{x\}$ is trivial.
\end{enumerate}
\end{fct}

Let $\frkM(X,G)$ be the moduli stack of $G$-torsors 
on a smooth projective curve $X$.
We have a natural morphism
$\Gr(X,G)_{x} \to \frkM(X,G)$
by forgetting the trivialization $\varphi$.
If $G$ is semisimple, then
Fact \ref{fct:DS} implies the following 
adelic description of $\frkM(X,G)$.
\[
 \frkM(X,G) \simeq 
 G(\shK_x)_{\text{out}}\backslash G(\shK_x)/G(\shO_x).
\]
Here 
$G(\shK_x)_{\text{out}}$
denotes the space of regular functions 
$X\setminus \{x\} \to G$,
which is naturally a subgroup of $G(\shK_x)$.

\subsection{Beilinson-Drinfeld Grassmannian}
\label{subsec:BDG}

Let $X$ be a smooth algebraic curve as before.
Denote by $\catSch$ the category of schemes over $\bbC$.

For $I \in \catS$,
consider the functor which maps $S \in \catSch$ to the data
$(f^I,\shP,\varphi)$ consisting of 
\begin{itemize}
\item 
a morphism $f^I:S \to X^I$ of schemes
\item
a $G$-torsor $\shP$ on $S \times X$
\item
a trivialization 
$\varphi$ of $\shP$ on 
$S \times X \setminus \bigl\{\Gamma(f^I_{i})\bigr\}_{i \in I}$,
where 
$f^I_i:S \to X$ is the composition of $f^I$ 
with the $i$-th factor projection $X^I \to X$,
and $\Gamma(s) \subset S \times X$ is the graph scheme of morphism 
$s: S \to X$.
\end{itemize}
a finite set $\{x_i\}_{i \in I}$ of points on $X$.

By Fact \ref{fct:DS},
this functor can be represented by an ind-scheme $\Gr(X,G)_{I}$.
For $S=\Spec k$, 
$\Gr(X,G)_{I}(S)$ is identified with 
\begin{align*}
\left\{\bigl(\shP, \{x_i\}_{i \in I}, \varphi\bigr)
       \, \big| \, 
       \shP \in \catM_G(X), \, x_i \in X, \, 
       \varphi: \text{trivialization of }
       \left.\shP\right|_{X\setminus \{x_i\}_{i \in I}} \right\}.
\end{align*}
We have a natural morphism
\[
 r^{(I)}: \, \Gr(X,G)_{I} \longto X^I,
\]
and the space $\Gr(X,G)_{x}$ in \eqref{eq:BDG:Gr-XGx} 
is the fiber of $r^{(I)}$ with $I=\{1\}$.

\begin{fct}\label{fct:BD}
The collection 
\[
 \Gr(X,G) \seteq \{\Gr(X,G)_{I}\}_{I \in \catS}
\]
has a structure of factorization monoid on $X$,
called the \emph{Beilinson-Drinfeld Grassmannian}.
\end{fct}

Here we explain the factorization isomorphism
only for $I = \{1,2\}$,
Denote by $\Gr_{\{x_i\}_{i \in I}}$ the fiber of $r^{(I)}$ 
over $\{x_i\}_{i \in I} \subset X^I$.
By Fact \ref{fct:DS}, in the case $|I|=1$, we have 
\[
 \Gr_{\{x\}} \simeq G(\shK)/G(\shO).
\]
If $x_1 \neq x_2$, then $\{x_1,x_2\} \in X^2 \setminus \Delta$,
and we have a morphism
$\Gr_{\{x_1,x_2\}} \to \Gr_{\{x_1\}} \times \Gr_{\{x_2\}}$ 
by restricting the data to $X \setminus x_i$.
We also have the other direction map 
$\Gr_{\{x_1\}} \times \Gr_{\{x_2\}} \to \Gr_{\{x_1,x_2\}}$
by gluing the $G$-bundles over $X\setminus \{x_1,x_2\}$. 
If $x_1=x_2$, then $\{x_1,x_2\} \in \Delta$,
and we have an identification
$\Gr_{\{x_1,x_2\}} \simto \Gr_{\{x_1\}}$.

The unit on $\Gr(X,G)$ is given by the trivial $G$-bundles.
Namely, define $u^{(I)}: X^I \to \Gr_{X,G,I}$ by
setting $u^{(I)}\bigl(\{x_i\}_{i \in I}\bigr)$ to 
be the trivial $G$-bundle with 
the obvious trivialization away from $x_i$'s.

$\Gr(X,G)$ has a natural factorization line bundle.
Assume the Lie algebra $\Lie(G)$ 
is equipped with an invariant inner product.
It induces a central extension
\[
 1 \longto \bbG_m \longto \wh{G} \longto G(\shK) \longto 1
\]
of algebraic group ind-schemes,
and hence a $\bbG_m$-torsor
\[
 \wh{G}(\shK)/G(\shO) \longto G(\shK)/G(\shO)
\]
over the affine Grassmannian.
It defines a factorization line bundle $\shL(G)$ over $\Gr(X,G)$.

\section{Linearization}
\label{sect:linear}

\subsection{Recollection on $\shD$-modules}
\label{subsec:ca-fa:Dmod}

Let us collect some notations on $\shD$-modules. 

Let $Z$ be a smooth algebraic variety over $\bbC$. 
%
Recall that the canonical sheaf 
\[
 \omega_Z \seteq \wedge^{\dim Z}\Omega_Z
\] 
has a natural right $\shD$-module structure determined by
\[
 \nu \cdot \tau := -\LieD_{\tau}(\nu),\quad
 \nu \in \omega_Z, \ \tau \in \Theta_Z \subset \shD_Z
\]
where $\LieD_{\tau}$ is the Lie derivative with respect to $\tau$.
Then the categories of left and right $\shD$-modules are equivalent 
under the functor
\[
 \shL \longmapsto \shL^r \seteq \omega_Z \otimes_{\shO_Z} \shL,
\] 
where the structure of right $\shD$-module on $\shL^r$ is determined by 
\[
 (\nu \otimes l)\cdot \tau \seteq
  (\nu \cdot \tau) \otimes l - \nu \otimes (\tau \cdot l).
\]
The inverse functor is given by
\begin{equation}\label{eq:d-mod:ell}
 \shM \longmapsto 
 \shM^{\ell} \seteq \omega_Z^{-1} \otimes_{\shO_Z} \shM.
\end{equation}

Hereafter the term a ``$\shD$-module" means a right $\shD$-module.

\subsection{Linearization of factorization spaces}
\label{subsec:fs:linearization}

Let $\shG$ be a factorization monoid over a smooth variety $X$.
Recall the morphisms $r^{(I)}: \shG_{I} \to X$ and 
$u^{(I)}: X \to \shG_{I}$.
Now consider the $\shO$-module 
\[
  \shA_{\shG,I} \seteq r^{(I)}_* u^{(I)}_{!} \omega_{X^I}.
\]
This sheaf can be considered as the space of delta functions on $\shG_I$ 
along the section $u^{(I)}(X^I)$.

The connection on $\shG$ along $X$ given in Fact \ref{fct:conn}
defines a right $\shD$-module structure on $\shA_{\shG,I}$
and the section $u^{(I)}$ defines 
an embedding $\omega_{X^I} \inj \shA_{\shG,I}$.
Then the axiom of factorization space implies

\begin{fct}\label{fct:fs-fa}
If $X$ is a smooth curve, then
the collection
$\{\shA_{\shG,I}^{\ell}\}_{I \in \catS}$ 
has a structure of factorization algebra on $X$,
and hence
$\shA_{\shG,\{1\}}$ 
has a structure of chiral algebra on $X$.
\end{fct}

Here we used the notation \eqref{eq:d-mod:ell}.
We also omit the definition of chiral and factorization algebras.
See \cite[\S\S3.1--3.4]{BD}.
We call the obtained chiral algebra $\shA_{\shG,\{1\}}$
the \emph{chiral algebra associated to $\shG$}.

Let us explain the twisted version of this construction.
Assume $X$ is a curve and 
a factorization space $\shG$ over $X$ 
is equipped with a factorization line bundle $\shL=\{\shL_I\}$.
Then consider the collection of $\shD$-modules 
\[
 \shA_{\shG,I}^{\shL} \seteq 
 r^{(I)}_* \bigl(\shL_I \otimes_{\shG_I} u^{(I)}_{!} ( \omega_{X^I})\bigr).
\]
The construction in Fact \ref{fct:fs-fa} can be applied to 
this family of $\shL_I$-twisted sheaves.
The result is 

\begin{fct}\label{fct:fs-fa:twist}
Assume that 
$\shG$ is a monoid over a smooth curve $X$,
and that a factorization linear bundle $\shL$ over $\shG$ is given.
Then $\shA^{\shL}_{\shG,\{1\}}$ 
has a structure of chiral algebra on $X$.
\end{fct}

We call the resulting chiral algebra
the \emph{$\shL$-twisted chiral algebra associated to $\shG$}.

\section{Lattice chiral algebra}
\label{sect:lattice}

Now we recall the factorization monoid
arising from the Picard functor
whose associated chiral algebra is the lattice chiral algebra
\cite[\S3.10]{BD}, \cite[Chap.\ 6]{G}, \cite[\S\S20.4.4--9]{FB}.

Following the symbols in \cite[Chap.\ 5]{FB},
let $V_{\bbZ} = \oplus_{n\in\bbZ} \pi_n$ be 
the lattice vertex algebra
associated to the lattice $\bbZ$.
Each component $\pi_n$ is the Fock representation of 
the Heisenberg algebra
generated by $\{a_n \mid n \in \bbZ\}$ and $1$ 
with the defining relation $[a_m,a_n]=m\delta_{m+n,0}$.
Recall that $V_{\bbZ}$ is actually a vertex \emph{super} algebra,

There is a universal chiral algebra $\shA_{\bbZ}$ 
in the meaning of \cite[\S3.3.14]{BD} corresponding to 
the lattice vertex algebra $V_{\bbZ}$.
We call it the \emph{lattice chiral algebra}.


\subsection{Some super language}

Below we will use some super language,
so let us collect the notation.

For a commutative ring $R$,
a super $R$-module means a $\bbZ/2\bbZ$-graded $R$-module.
Denote by $p$ the parity of homogeneous elements of $M$.
Consider the tensor product of super $R$-modules 
with the commutativity constraint 
$a \otimes b =(-1)^{p(a) p(b)} b \otimes a$.
A super $R$-line is an invertible object of 
the monoidal category formed by super $R$-modules with the tensor product 
described as above.
Super $R$-lines constitute a Picard groupoid.

One can consider a similar setting for 
any $\bbC$-linear category $\catM$.
Namely a $\bbZ/2\bbZ$-graded object form a $\bbC$-category $\catM^s$,
and if $\catM$ is a monoidal category, then
$\catM^s$ is also a super monoidal category with Koszul rule of sings
as described above.

Below we will use super line bundles on a scheme (or stack) $X$.
It is nothing but a super object in 
the monoidal category $\Pic(X)$.
In particular, replacing line bundles by super line bundles in 
Definition \ref{dfn:line},
we have the notion of \emph{factorization super line bundle} 
on a factorization space.

\subsection{Abelian Beilinson-Drinfeld Grassmannian and lattice chiral algebras}

Consider the Beilinson-Drinfeld Grassmannian $\Gr(X,G)$ 
with $G=\GL(1)$.
Since $\GL(1)$-torsor is nothing but a line bundle,
we have $\frkM(X,\GL(1)) = \Pic(X)$.
For $I \in \catS$, we also have 
\[
 \Gr(X,\GL(1))_{I} =
 \left\{(\shL, S, \varphi) \mid 
  \shL \in \Pic(X), \, S \in X^I, \, 
  \varphi: \text{trivialization of } \shL|_{X \setminus S}
 \right\}.
\]
%
%
%
%
%
%


One can naturally construct a factorization \emph{super} line bundle 
on $\Gr(X,\GL(1))$ by the theta line bundle $\Theta$ on $\Pic(X)$.
Choosing a theta characteristic $\Omega^{1/2}_X$,
we can describe the fiber of $\Theta$ on $\shL \in \Pic(X)$ as 
\begin{equation}\label{eq:Theta}
 \left.\Theta \right|_{\shL} \dsimeq 
 \det H^*(\shL \otimes_{\shO_X} \Omega_X^{1/2}).
\end{equation}
Considered as an odd line bundle on $\Pic(X)$,
it induces a factorization super line bundle on $\Gr(X,\GL(1))$.
Denote it by the same letter $\Theta=\{\Theta_I\}$.

Recall that the chiral algebra associated to a factorization monoid
is constructed along the unit $u^{\{1\}}: X \to \shG_{\{1\}}$.
In the present situation $\shG=\Gr(X,\GL(1))$,
the unit is given by 
$x \mapsto (\shO_X,\varphi)$ with $\varphi$ the obvious trivialization 
away from $x$.
Now we modify this construction and consider
\[
 u(n)_x: \, X \longto \Gr(X,\GL(1))_{x}, \quad
 x \longmapsto (\shO_X(n x),\varphi)
\]
for $n \in \bbZ$.
It induces $u(n)^{\{1\}}: X \longto \Gr(X,\GL(1))_{\{1\}}$
and $u(n)^{(I)}: X \longto \Gr(X,\GL(1))_{I}$ for any $I \in \catS$.
Now define the $\shO$-module
\[
 \shA_{\bbZ,I} \seteq 
 \bigoplus_{n \in \bbZ} 
 r_*^{(I)}\left(\Theta_I \otimes u(n)^{(I)}_!(\omega_{X^I})\right).
\]
The factorization structure makes $\shA_{\bbZ,\{1\}}$ 
into a chiral algebra. 
Denote it by $\shA_{\bbZ}$.

\begin{fct} 
\label{fct:lattice}
The chiral algebra $\shA_{\bbZ}$ is isomorphic to the 
lattice chiral algebra.
\end{fct}

For later use let us restate this result.
Set 
\[
 \shG(X,\bbZ)_I \seteq 
 \{(\shL \otimes \shO_X(*S),S,\varphi) \mid 
   \shL \in \Pic(X), S \in X^I,
   \varphi: \text{trivialization of } \shL|_{X \setminus S}
 \}
\]
Here $\shO_X(*S)$ means $\shO_X(n S)$ for some $n \in \bbZ$.
Then
$\shG(X,\bbZ)$ is a factorization monoid with the unit 
$u:X \to \shG(X,\bbZ)_I$ given by the trivial line bundle.
$\Theta$ is a factorization super line bundle on $\shG(X,\bbZ)$.
Now $\shA_{\bbZ}$ is isomorphic to the $\Theta$-twisted 
chiral algebra associated to $\shG(X,\bbZ)$.
Thus Fact \ref{fct:lattice} can be restated as 

\begin{prop}\label{prop:lattice}
The $\Theta$-twisted 
chiral algebra associated to $\shG(X,\bbZ)$
is isomorphic to the lattice chiral algebra.
\end{prop}

\begin{rmk*}
The $\Theta$-twisted chiral algebra associated to $\Gr(X,\GL(1))$ 
is nothing but the Heisenberg chiral algebra $\pi^{\ch}_0$,
and the corresponding vertex algebra $\pi_0$ is the $0$-th part 
of the lattice vertex algebra $V_{\bbZ} = \oplus_{n \in \bbZ} V_n$. 
\end{rmk*}

\section{Clifford chiral algebras}
\label{sect:cl}

\subsection{Clifford bundles}

Here we collect some notations on the Clifford algebra.
See \cite[\S IX.9]{B} for the detail.

Let us denote by $\Cl(Q)$ 
the \emph{Clifford algebra} associated to non-degenerate quadratic form $Q$ 
on a linear space $V$ of finite dimension.
Denote by $\Cl^{+}(Q)$ the even part and $\Cl^{-}(Q)$ 
the odd part of $\Cl(Q)$.
The involution $u \mapsto \pm u$ for $u \in \Cl^{\pm}(Q)$ 
is denoted by $\alpha$.
The \emph{Clifford group} $\C(Q)$ and 
the \emph{special Clifford group} $\SC(Q)$ are defined to be
\[
 \C(Q) \seteq 
 \{u \in \Cl^{\times}(Q) \mid \alpha(u) V u^{-1} \subset V\},
 \quad
 \SC(Q) := \C(Q) \cap \Cl^{+}(Q),
\]
where $\Cl^{\times}(Q)$ is the group of units in $\Cl(Q)$.
$\SC(Q)$ is a connected reductive group.
The \emph{spinor norm} is denoted by
$N: \SC(Q) \to \bbC^{\times}$.

Using the notations in \S\ref{subsec:Bun_G},
we have a moduli stack $\frkM(\SC(Q),X)$ 
of $\SC(Q)$-torsors on a smooth algebraic curve $X$.
Applying Fact \ref{fct:BD} to $G=\SC(Q)$.
we have a factorization monoid $\Gr\left(X,\SC(Q)\right)$.
The spinor norm $N$ on $\SC(Q)$ induces 
a factorization super line bundle 
$\shL(\SC(Q))$ by a similar argument 
as in the latter part of \S\ref{subsec:BDG}.


\subsection{Clifford chiral algebras}

Let us introduce the Clifford chiral algebra
following \cite[\S3.8.6]{BD} with minor modifications.

Let $\shV, \shV'$ be vector $\shD_X$-bundles 
(locally projective right $\shD_X$-module of finite rank) 
and assume that we have a skew-symmetric pairing 
$\pair \in P^*_2(\{\shV,\shV'\},\omega_X)$
in the $*$-pseudo-tensor category \cite[Chap.\ 2]{BD}
such that it vanishes on $\shV$ and $\shV'$ respectively.
Then $\shV \oplus \shV' \oplus \omega_X$ 
is a $\Lie^*$ algebra 
with the commutator $\pair$.
The resulting chiral algebra is denoted by $\Wch(\shV,\shV',\pair)$.

For a positive integer $n$,
let $\shV_o :=  \shO_X^{\oplus n}$ considered as $\shO$-module,
and $\shV := \shV_o \otimes_{\shO_X} \shD_X$ considered as $\shD$-module.
Set $\shV^{\circ} := \shHom^*(\shV,\shO_X)$ with $*$ denoting the linear dual.
As in \cite[\S2.2.16]{BD}, 
the linear pairing induces a non-degenerate pairing 
$\pair \in P^*_2(\{\shV[1],\shV^{\circ}[-1]\},\omega_X)$ 
with $[\pm1]$ the shift as complexes.
Now we set 
\[
 \Clch(n) \seteq \Wch(\shV[1],\shV^{\circ}[-1],\pair)
\]
and call it the \emph{fermionic chiral algebra}.
Thus, strictly speaking, it is a chiral \emph{super} algebra and
generators behave as odd elements.
It is a universal chiral algebra 
associated to the free fermionic vertex algebra $\Wedge$.


\subsection{}

Let us set $\SC(m) := \SC(Q)$ 
with $Q$ the standard quadratic form $z_1^2 + \cdots +z_m^2$ 
on $\bbC^{m}$. 
Then we have the factorization monoid $\Gr(X,\SC(m))$ 
with the factorization super line bundle $\shL(\SC(m))$.
Let us denote by $\shA(\SC(2m))$
the $\shL(\SC(m))$-twisted chiral algebra 
associated to $\Gr(X,\SC(m))$.

Now we have

\begin{prop}\label{prop:cl}
$\shA(\SC(2))$ is isomorphic to $\Clch(1)$.
\end{prop}

The proof is a variant of the argument in \cite[Chap.\ 20]{FB}.
It says that for a semi-simple $G$,
the $\shL(G)^{\otimes k}$-twisted chiral algebra 
associated to the Beilinson-Drinfeld Grassmannian $\Gr(X,G)$ 
is the affine chiral algebra $\shV^{\ch}_k(\g)$ of level $\bbC$.
Here $\shV^{\ch}_k(\g)$ is the chiral algebra counterpart of 
the affine vertex algebra $V_k(\g)$ of level $k$
with $\g := \Lie(G)$,
and $\shL(G)$ is the factorization line bundle 
induced by the normalized invariant inner product
(see the latter part of \S\ref{subsec:BDG}).

\begin{proof}
Denote by $\Wedge$ the free fermionic vertex super algebra
(following the symbol in \cite[\S5.3]{FB}).
Thus as a linear space $\Wedge$ is the fermionic Fock space 
of the infinite-dimensional Clifford algebra 
$\Cl_\infty$.
The algebra $\Cl_\infty$ is 
generated by $\{\psi_n, \psi^*_n  \mid n \in \bbZ\}$
with the relation
\[
 [\psi_m,\psi_n]_+  =  [\psi^*_m,\psi^*_n]_+ = 0,\quad
 [\psi_m,\psi^*_n]_+  = \delta_{m+n,0}.
\] 
The $\Cl_{\infty}$-module $\Wedge$ is generated 
by the element $\ket{0}$ such that  
\[
 \psi_n \ket{0} = 0 \ (n \ge 0),\quad
 \psi^*_n \ket{0} = 0 \ (n > 0).
\]

By the description of linearization procedure,
we have an identification
$\shA(\SC(2))_x \simeq \Cl_{\infty} \otimes_{\Cl_\infty^+} \bbC$ 
for any $x \in X$.
Here $\bbC$ is the one-dimensional representation of $\Cl^+_{\infty}$
corresponding to the fiber of $\shL(\SC(2))_{\{1\}}$.
Thus $\shA(\SC(2))_x$ is identified with 
the free fermion vertex algebra $\Wedge$ as linear space.

The factorization structure induces a connection on 
$\shA^{\ell}(\SC(2))$ giving a left $\shD$-module structure,
Denote by $\psi(z) z^{1/2} d z$ and $\psi^*(z) z^{-1/2} d z$  
the generating local sections of $\shA^{\ell}(\SC(2))$.
Here $z^{\pm 1/2}$ corresponds to 
the theta characteristic appearing in \eqref{eq:Theta}.
$\psi(z)$ and $\psi^*(z)$ are considered to be odd fields.
Then as in the affine chiral algebra case,
we have $(z-w)[\psi_i(z),\psi_j^*(w)]_+ =0$,
which is nothing but 
the Clifford chiral super algebra relation.
\end{proof}

\section{Boson-Fermion correspondence via factorization spaces}
\label{sect:main}

\begin{thm}\label{thm:main}
We have an isomorphism of factorization monoids
\[
 \shG(X,\bbZ) \dsimeq \Gr(X,\SC(2)).
\]
Under this isomorphism we have an identification
of factorization super line bundles 
\[
 \Theta \dsimeq \shL(\SC(2)).
\]
\end{thm}

Then Propositions \ref{prop:lattice} and \ref{prop:cl} 
yield $\shA_{\bbZ} \simeq \Clch(1)$ as chiral algebras.
In terms of  vertex algebra, it means the isomorphism
$V_{\bbZ} \simeq  \Wedge$
which is nothing but the boson-fermion correspondence.

\begin{proof}[{Proof of Theorem \ref{thm:main}}]
Recall the notations of $G$-bundles in \S\ref{subsec:Bun_G}.
Since we have
$\SC(2) \simeq \GL(1)^{\times 2}$,
it holds that
\[
 \frkM(\SC(2),X) \dsimeq 
 \frkM(\GL(1),X)^{\times 2} \dsimeq 
 \Pic(X)^{\times 2}.
\]
Then by Fact \ref{fct:DS} on the trivialization we have
\[
 \Gr(X,\SC(2))_{x} \dsimeq 
 \bigl(\Gr(X,\GL(1))_{x}\bigr)^{\times 2} \dsimeq 
 \shG(X,\bbZ)_{x}
\]
for each $x \in X$.
It yields the isomorphism
$\Gr(X,\SC(2))_{I} \simeq \shG(X,\bbZ)_{I}$
of ind-schemes for any $I \in \catS$.

The factorization structures on $\Gr(X,\SC(2))$ 
and $\shG(X,\bbZ)$ are the same essentially because 
they are both Beilinson-Drinfeld Grassmannians.
The units coincide because they are given by trivial bundles.
Thus we have $\shG(X,\bbZ) \simeq \Gr(X,\SC(2))$.

Recall that the factorization super line bundle $\shL(\SC(2))$ 
comes from the spinor norm $N$ on $\SC(2)$.
Under the isomorphism $\SC(2) \simeq \GL(1)^{\times 2}$,
$N$ is given by $N(x,y) = x y$,
and it can be identified with the linear pairing.
Now the isomorphism $\Theta \simeq \shL(\SC(2))$ 
is an obvious one.
\end{proof}

\subsection*{Acknowledgements.} 
The author is supported by the Grant-in-aid for 
Scientific Research (No.\ 16K17570), JSPS.
This work is also supported by the 
JSPS for Advancing Strategic International Networks to 
Accelerate the Circulation of Talented Researchers
``Mathematical Science of Symmetry, Topology and Moduli, 
  Evolution of International Research Network based on OCAMI''.

\end{document}